\documentclass[12pt,reqno]{amsart}
\usepackage{amssymb}
\usepackage{a4wide}
\usepackage{mathrsfs}
\def\mathcal{\mathscr}
\usepackage{newsymbol}
\usepackage{verbatim}
\usepackage[active]{srcltx}
\usepackage[colorlinks,linkcolor={blue},
citecolor={blue},urlcolor={red},]{hyperref}
\usepackage{hyperref}
\let\emptyset \undefined
\let\ge       \undefined
\let\le       \undefined
\newsymbol\le          1336  \let\leq\le
\newsymbol\ge          133E  \let\geq\ge
\newsymbol\emptyset    203F
\newsymbol\notle       230A
\newsymbol\notge       230B
\theoremstyle{plain}
\newtheorem{theorem}{Theorem}[section]
\theoremstyle{plain}
\newtheorem{corollary}[theorem]{Corollary}

\newtheorem{proposition}[theorem]{Proposition}

\newtheorem{hypothesis}[theorem]{Hypothesis}
\theoremstyle{remark}
\newtheorem{remark}[theorem]{Remark}

\numberwithin{equation}{section}
\def\R{{\mathbb R}\,}

\def\E{{\mathbb E}\,}
\def\P{{\mathbb P}\,}
\def\la{\left(}
\def\ra{\right)}

\begin{document}
\title[]{Exponential ergodicity of semilinear equations driven by L\'evy processes in Hilbert spaces}
\author{A. Chojnowska-Michalik}
\address{Faculty of Mathematics and Computer Science\\
University of \L \'od\'z\\
Banacha 22, 90238 \L \'od\'z\\
 Poland}
\email{anmich@math.uni.lodz.pl}
\author{B. Goldys}
\address{School of Mathematics and Statistics\\
The University of Sydney\\
Sydney 2006, Australia}
\email{beniamin.goldys@sydney.edu.au}
\thanks{This work was partially supported by the ARC grant DP120101886 and by the Polish Ministry of Science and Higher Education grant "Stochastic
equations in infinite dimensional spaces" N N201 419039}
\begin{abstract}
We study convergence to the invariant measure for a class of semilinear
stochastic evolution equations driven by L\'evy noise, including the case of cylindrical
noise. For a certain class of equations we prove the exponential rate of
convergence in the norm of total variation. Our general result is applied to a number of specific equations driven by cylindrical symmetric $\alpha$-stable noise and/or cylindrical Wiener noise. We also consider the case of a "singular" Wiener process with unbounded covariance operator. In particular, in the equation with 
diagonal pure $\alpha$-stable cylindrical noise introduced by Priola and Zabczyk we generalize results in \cite{PSXZ}. In the proof  we use an idea of Maslowski and Seidler from \cite{MS}.   
\end{abstract}
\maketitle
\tableofcontents
\section{Introduction}
This work is dedicated to Prof. Jerzy Zabczyk, our Master and Teacher. 
\par\bigskip\noindent
Let $H$ be a real separable Hilbert space. In this paper we are concerned with the ergodic behaviour of the semilinear stochastic equation   
\begin{equation}\label{1}
\begin{cases} dX_t = \la A X_t +F\la X_t\ra\ra dt + dZ_t, &\quad t>0,\\
X_0 =x\in H, &{}\end{cases}
\end{equation}
driven by, possibly cylindrical, L\'evy noise $Z$. We assume that $A$ is a generator of a $C_0$-semigroup $(S_t)_{t\geq 0}$ of bounded linear operators on $H$ and $F:H\to H$ is a Lipschitz mapping. Under the conditions specified in Section \ref{main_result} equation \eqref{1} has a unique solution $\la X_t^x\ra$ for every $x\in H$. Our main question is to provide conditions for the exponential convergence of transition measures $\mu_t^x(\cdot)=\P\la X_t^x\in\cdot\ra$ to a unique invariant measure $\mu$ in the norm of total variation. This question has been thoroughly investigated in the case when $Z=BW$ with a certain cylindrical Wiener process $W$ and $B$ a linear operator, see for example \cite{GM06,GM,m93}. In those works more general, non-Lipschitz drifts $F$ were permitted. The two main ingredients for the proofs were the strong Feller property and topological irreducibility of the transition measures $\mu_t^x$. \\ 
\par\noindent
In this paper we extend the aforementioned approach to equation \eqref{1} with general (subject to some natural conditions) L\'evy noise $Z$ and bounded Lipschitz drift $F$. Exponential ergodicity for equation \eqref{1} was investigated recently in \cite{PSXZ,PXZ}.  In Theorem \ref{th1} of Section \ref{main_result} we show that if the semigroup $\la S_t\ra$ is exponentially stable and compact and the measures $\mu_t^x$ are  strong Feller and irreducible then the convergence to the invariant measure holds in the total variation norm and moreover the rate of convergence is exponential. The proof is based on an idea developed in \cite{MS}. In Section \ref{examples} we discuss three examples. In subsection \ref{stable} we apply our general result to equation (\ref{1}) with the noise $Z$ being a symmetric $\alpha$-stable cylindrical noise. In Corollary \ref{cor1} we show that Theorem \ref{th1} yields a strengthening of the similar result obtained in \cite{PSXZ} and of an earlier result in \cite{PXZ}. In the proof we use the fundamental results of Priola and Zabczyk proved in \cite{PrZ3}. In subsection \ref{sec-gaussian} we consider L\'evy noise with both, the jump part and the Gaussian part, non-vanishing. Finally, in Section \ref{singular} we consider the stochastic heat equation driven by a mixture of a cylindrical symmetric $\alpha$-stable noise and a singular Wiener noise with unbounded covariance operator. 
\par\medskip\noindent
\textbf{Notations. }In what follows, the norm in the Hilbert space $H$ is denoted by $|\cdot|$. The norm in any other Banach space $E$ will be denoted as $|\cdot|_E$. 
\par
The Hilbert-Schmidt norm of an operator $K:H\to H$ will be denoted as $\|K\|_{HS}$. 
\section{Formulation of the problem and main result}\label{main_result}
We start with the set of conditions that are assumed to hold throughout the paper. We assume that $Z= (Z_t)_{t\geq 0}$ is a L\'evy process defined on a stochastic basis $(\Omega, \mathcal  F, (\mathcal  F_t)_{t\geq 0}, \mathbb P)$ satisfying the usual conditions and with values in a Hilbert space $U$, such that $H\subset U$ with continuous and dense embedding. 
More precisely, we assume that 
$$Z_t =ta +Q^{1/2}W_t +Y_t,\quad t\geq 0,$$
where $a\in U$ is a fixed vector, $W$ is a standard cylindrical Wiener process in $H$, the process $Y$ is a $U$-valued, pure jump L\'evy process, independent of $W$, and $Q:\mathrm{dom}(Q)\subset H\to H$ is a selfadjoint and nonnegative operator.\\
We will assume that the semigroup $\la S_t\ra$ enjoys the following properties. 
\begin{hypothesis}\label{h1}
\begin{enumerate}
\item
For every $t\in(0,1]$, the operator $S_t$ has an extension (still denoted by $S_t$) to a bounded operator $S_t:U\to H$. 
\item
We have \begin{equation}\label{eqh1}
\int_0^1\left|S_sa\right|\,ds<\infty.
\end{equation}
For each $t>0$ the operator $S_tQ^{1/2}$ extends to a bounded operator on $H$ and 
\begin{equation}\label{eqh2}
\int_0^1\left\|S_sQ^{1/2}\right\|^2_{HS}\,ds<\infty. 
\end{equation}
\end{enumerate}
\end{hypothesis}
\begin{hypothesis}\label{h2}
For every $t\in(0,1]$ the stochastic integral 
\[\int_0^tS_sdY_s,\]
is a well defined $H$-valued random variable defined as a limit in probability of the corresponding Stieltjes sums (see e.g. \cite{Ch} and \cite{PZ}). Moreover, there exists $p>0$ such that 
\begin{equation}\label{eqh3}
\gamma_p=\sup_{t\in(0,1]}\E\left|\int_0^tS_sdY_s\right|^p<\infty,
\end{equation} 
and 
\begin{equation}\label{eqh4}
\lim_{t\downarrow 0}\E\left|\int_0^tS_sdY_s\right|^p=0.
\end{equation}
\end{hypothesis}
The remark below provides a convenient sufficient condition for Hypothesis \ref{h2} to hold. It is proved in Section \ref{proofs}. 
\begin{remark}\label{22a}
If $U=H$ and $\E\left|Y_1\right|^p<\infty$ for a certain $p>0$ then Hypothesis \ref{h2} is satisfied with this $p$. 
\end{remark}
Finally, we impose the standard Lipschitz condition on the nonlinear term $F$.
\begin{hypothesis}\label{h3}
The function $F:H\to H$ is Lipschitz continuous and bounded. 
\end{hypothesis}
\begin{proposition}\label{pr1}
Assume Hypothesis \ref{h2} and let 
\[Y_A(t)=\int_0^tS_{t-s}dY_s,\quad t\ge 0.\]
Then the following hold. 
\begin{enumerate}
\item
For any $t\ge 0$ there exists a constant $C\la t,p,A\ra$ such that 
\begin{equation}\label{eq2}
\E\left|Y_A(t)\right|^p\le\gamma_pC(t,p,A). 
\end{equation}
\item
The process $Y_A$ is continuous in $p$-th mean.
\item
The process $Y_A$ admits a predictable modification.
\end{enumerate}
\end{proposition}
\begin{remark}\label{rem1}
Clearly, Proposition \ref{pr1} holds true for a more more general process 
\[\begin{aligned}
Z_A(t)&=\int_0^t S_{t-s} dZ_s\\
&= \int_0^t S_s a\,ds + \int_0^t S_{t-s} Q^{1/2} dW_s+ \int_0^t S_{t-s}d Y_s, \quad t\geq 0,
\end{aligned}\]
with 
\[\tilde{\gamma}_p=\sup_{t\in(0,1]}\E\left|\int_0^tS_sdZ_s\right|^p<\infty,\]
by \eqref{eqh1}, \eqref{eqh2} and \eqref{eqh3}. We will always consider a predictable modification of $Z_A$. 
\end{remark}
Let us recall that the process $X^x$ is a mild solution to equation \eqref{1} if it is $H$-valued, predictable and 
\begin{equation}\label{1a}
X^x_t=S_tx+\int_0^tS_{t-s}F\la X_s^x\ra ds+Z_A(t),\quad t\ge 0,\quad\P-a.s.
\end{equation}
We will use the notation 
\[\mu_t^x\quad\mathrm{for\,\, the\,\,law\,\,of\,\,the\,\,random\,\,variable}\quad X_t^x\in H.\]
\begin{proposition}\label{pr2}
Assume Hypotheses \ref{h1}-\ref{h3}. Then for every $x\in H$ equation \eqref{1} has a unique mild solution $X^x$. The process $X^x$ is Markov and Feller in $H$, 
\[\E\left|X_t^x\right|^p<\infty,\quad t\ge 0,\]
and the process $X^x$ is $t$-continuous in $p$-th mean. 
\end{proposition}
Now we can formulate the main result of this paper. In the theorem below we denote by $\|\nu\|_{var}$ the total variation norm of a signed Borel measure $\nu$ on $H$. For the definitions of the strong Feller property and open set irreducibility the reader may consult \cite{DZ1}. 
\begin{theorem}\label{th1}
Assume Hypotheses \ref{h1}-\ref{h3}. Moreover assume that the semigroup $\la S_t\ra$ is  compact and exponentially stable. Assume also that the solution $X$ to \eqref{1} is strongly Feller and open set irreducible. Then there exists a unique invariant measure $\mu$ for \eqref{1} and there exist constants $\beta>0$ and $C>0$ such that 
\begin{equation}\label{eq3}
\left\|\mu_t^x-\mu\right\|_{\mathrm{var}}\le Ce^{-\beta t}\la 1+|x|^p\ra,\quad x\in H,\,t\ge 0.
\end{equation}
\end{theorem}
\begin{remark}\label{rem2}
Let $X_t^\nu$ be the solution to \eqref{1} starting from an $\mathcal F_0$-measurable random variable $X_0$ such that $\mathrm{Law}\la X_0\ra=\nu$. Then inequality \eqref{eq3} extends  by standard arguments to a more general one
\begin{equation}\label{eq3a}
\left\|\mu_t^\nu-\mu\right\|_{\mathrm{var}}\le Ce^{-\beta t}\la 1+\int_H|x|_H^p\nu(dx)\ra,\quad t\ge 0,
\end{equation}
where $\mu_t^\nu$ is the law of the random variable $X_t^\nu$. 
\end{remark}
\section{Examples}\label{examples}
\subsection{Model with cylindrical $\alpha$-stable noise}\label{stable}
In this section we apply Theorem \ref{th1} to equation \eqref{1} driven by a cylindrical symmetric $\alpha$-stable noise $Z$ with $\alpha\in(0,2)$, introduced by Priola and Zabczyk in \cite{PrZ3} and later studied in \cite{PSXZ} and \cite{PXZ}. 
\par
Let $\mathrm{dim} (H)=\infty$ and let $\left\{e_k;\, k\ge 1\right\}$ be a CONS in $H$. We will  identify $H$ with the space $l^2$ of square-summable sequences. Let $\la\rho_k\ra$ be an arbitrary sequence of positive real numbers. Then we can define the Hilbert space $l^2_\rho$ of all sequences $x=\la x_k\ra$ such that 
\[|x|^2_\rho=\sum_{k=1}^\infty x_k^2\rho_k^2<\infty,\]
and we will denote by $\langle\cdot,\cdot\rangle_\rho$ the corresponding inner product.\\
The equation 
\begin{equation}\label{52}
\left\{\begin{aligned}
dX_t&=\la AX_t+F\la X_t\ra\ra dt+dZ_t,\,\, t>0,\\
X_0&=x\in H,
\end{aligned}\right.
\end{equation}
is studied in \cite{PrZ3} and in \cite{PXZ} in the space $H=l^2$ under the set of assumptions listed below.
\begin{hypothesis}\label{hpz}
\begin{enumerate}
\item
Let $A$ be a selfadjoint operator $A:\mathrm{dom}(A)\subset H\to H$ such that 
\[Ae_k=-\lambda_k e_k,\quad k\ge 1,\]
with 
\[0<\lambda_1\le\lambda_2\le\cdots,\quad \lambda_k\to\infty.\]
\item
For $\alpha\in(0,2)$ we define a cylindrical $\alpha$-stable process 
\[Z(t)=\sum_{k=1}^\infty b_kZ_k(t)e_k,\quad t\ge 0,\]
where $\la Z_k\ra$ is a sequence of real-valued symmetric and independent $\alpha$-stable processes and $\la b_k\ra$ is a sequence of positive reals such that 
\begin{equation}\label{53} 
\sum_{k=1}^\infty\frac{b_k^\alpha}{\lambda_k}<\infty.
\end{equation}
\item
The mapping $F:H\to H$ is Lipschitz continuous and bounded.
\item
There exist constants $\theta\in (0,1)$ and $\tilde{c}>0$ such that 
\[b_k\ge \tilde{c}\lambda_k^{-\theta+\frac{1}{\alpha}},\quad k\ge 1.\]
\end{enumerate}
\end{hypothesis}
Let for a fixed $\sigma\ge\frac{1}{\alpha}$, 
\[U_\sigma=l^2_{\rho(\sigma)},\quad \rho(\sigma)=\la\lambda_k^{-\sigma}\ra.\]
Clearly, $H$ is continuously and densely imbedded into $U_\sigma$.\\
We shall show that Hypothesis \ref{hpz} implies the assumptions of Proposition \ref{pr2} and of Theorem \ref{th1}.\\
First, note that $Z$ is an $U_\sigma$-valued L\'evy process. Indeed, by Proposition 3.3 in \cite{PrZ3}, $Z$ is a L\'evy process with values in $l_\rho$ with $\rho=\la\rho_k\ra$ such that 
\[\sum_{k=1}^\infty b_k^\alpha\rho_k^\alpha<\infty,\]
hence with values in $U_\sigma$ by \eqref{53}.\\
 Next, the operator $A$ with the domain $\mathrm{dom}(A)=l_\lambda^2$, where $\lambda=\la\lambda_k\ra$ generates on $H$ a symmetric compact semigroup $\la S_t\ra$ with 
\[S_te_k=e^{-\lambda_kt},\quad k\ge 1,\,t\ge 0.\]
In particular, 
\begin{equation}\label{exp}
\left\|S_t\right\|=e^{-\lambda_1t},\quad t\ge 0.
\end{equation}
Therefore, 
\begin{equation}\label{33a}
\mathrm{the\,\,additional\,\,assumptions\,\,on\,\,}\la S_t\ra\mathrm{\,\,in\,\,Theorem\,\,\ref{th1}\,\,hold.}
\end{equation} 
Fix $t>0$ and let $x=\la x_k\ra\in H\subset U_\sigma$. Then 
\begin{equation}\label{55}
\left|S_tx\right|^2_H=\sum_{k=1}^\infty e^{-2\lambda_kt}x_k^2=\sum_{k=1}^\infty\lambda_k^{2\sigma}e^{-2\lambda_kt}\frac{x_k^2}{\lambda_k^{2\sigma}}.
\end{equation}
Since $\lambda_k\nearrow\infty$,
\[c(t)=\sup_k\la\lambda_k^\sigma e^{-\lambda_kt}\ra<\infty\]
and \eqref{55} yields 
\[\left\|S_t\right\|_{U_\sigma\to H}\le c(t),\]
and the proof of part (1) of Hypothesis \ref{h1} is complete.\\ 
Since $a=0$ and $Q=0$ all the conditions of Hypothesis \ref{h1} are satisfied.\\
Next, if parts (1) and (2) of Hypothesis \ref{hpz} hold then it follows from Proposition 4.4 of \cite{PrZ3} that for any $t>0$ the stochastic integral $\int_0^tS_sdZ_s$ takes values in $H$, hence the first part of Hypothesis \ref{h2} is satisfied and additionally for any $p\in(0,\alpha)$ 
\begin{equation}\label{56}
\E\left|\int_0^tS_sdZ_s\right|^p\le\tilde{c_p}\la\sum_{k=1}^\infty b_k^\alpha\frac{1-e^{-\alpha\lambda_k t}}{\alpha\lambda_k}\ra^{p/\alpha},
\end{equation}
where the constant $\tilde{c_p}$ depends on $p$ only.  We note that the last estimate is an immediate consequence of Theorem 4.6 in \cite{PrZ3}. Note also that by \eqref{53} the RHS of \eqref{56} is finite uniformly in $t>0$, hence our condition \eqref{eqh3} is satisfied with 
\[\gamma_p\le \tilde{c_p}\la\frac{1}{\alpha}\ra^{p/\alpha}\la\sum_{k=1}^\infty\frac{b_k^\alpha}{\lambda_k}\ra^{p/\alpha}.\]
Moreover, using the Dominated Convergence we find that the RHS of \eqref{56} tends to 0 as $t\downarrow 0$ which yields our condition \eqref{eqh4}. Finally, all the conditions of Hypothesis \ref{h2} hold.\\
Hypothesis \ref{h3} is identical with part 3 of Hypothesis \ref{hpz}.\\
Thus we have proved that under assumptions (1)-(3) of Hypothesis \ref{hpz}, our Hypotheses \ref{h1}-\ref{h3} are satisfied and thereby as an immediate consequence of Proposition \ref{pr2} and Theorems 5.5 and 5.7 in \cite{PrZ3} we obtain 
\begin{corollary}\label{cor1}
Assume that Hypothesis \ref{hpz} holds. Then for every $x\in H$ there exists a unique mild solution $X^x$ to equation \eqref{52} which is an $H$-valued Markov process, irreducible and strong Feller. Moreover, for any $p\in(0,\alpha)$
\[\E\left|X_t^x\right|^p<\infty,\quad t\ge 0.\]
\end{corollary}
Finally, taking into account \eqref{33a} and \eqref{exp} we note that by Hypothesis \ref{hpz} all the assumptions of our Theorem \ref{th1} are satisfied and we obtain immediately the following strengthening of Theorem 2.8 in \cite{PSXZ} and of Theorems 2.2 and 2.3 in \cite{PXZ}.
\begin{corollary}\label{cor2}
Under Hypothesis \ref{hpz} there exists a unique invariant measure for \eqref{52}and for any $p\in(0,\alpha)$ there exist constants $\beta>0$ and $C>0$ such that 
\begin{equation}\label{35}
\left\|\mu_t^x-\mu\right\|_{var}\le Ce^{-\beta t}\la 1+|x|_H^p\ra,\quad t\ge 0,\, x\in H.
\end{equation}
\end{corollary}
\begin{remark}\label{34}
Ascertainment \eqref{35} is proved in Theorem 2.8 in \cite{PSXZ} under Hypothesis \ref{hpz} and an additional assumption that for some $\epsilon>0$ 
\begin{equation}\label{36}
\sum_{k=1}^\infty\frac{b_k^\alpha}{\lambda_k^{1-\alpha\epsilon}}<\infty.
\end{equation}
Let us note, that \eqref{36} ensures that the solution $X$ to \eqref{52} takes values in a linear subspace compactly imbedded into $H$, while \eqref{53} of Hypothesis \ref{hpz} is equivalent to the fact that $X$ is $H$-valued only (Comp. Remark 2.3 in \cite{PSXZ}). 
\par\noindent
Corollary \ref{cor2} strengthens also Theorems 2.2 and 2.3 from an earlier paper \cite{PXZ}.  In that paper additionally to Hypothesis \ref{hpz} the nonlinearity $F$ must be sufficiently small and only weak convergence to the invariant measure $\mu$ with exponential rate is obtained. 
\end{remark}
\subsection{Linear equation driven by L\'evy noise with Gaussian component}\label{sec-gaussian}
In this section we consider equation \eqref{1} with $F=0$, i.e. 
\begin{equation}\label{1b}
\begin{cases} dX_t = \la A X_t +a\ra dt+Q^{1/2}dW_t+ dY_t, &\quad t>0\\
X_0 =x\in H, &{}\end{cases}
\end{equation}
under Hypotheses \ref{h1} and \ref{h2}. By \eqref{eqh2} the operator $Q_t:H\to H$ defined as
\[Q_tx=\int_0^tS_sQS_s^\star x\,ds,\quad x\in H,\]
is nonnegative and trace class for every $t>0$. 
We will also assume the following 
\begin{hypothesis}\label{h3b1}
\[S_t(H)\subset Q_t^{1/2}(H),\quad t>0.\]
\end{hypothesis}
Let us recall that Hypothesis \ref{h3b1} holds if and only if the solution $X^{OU}$ (Gaussian Ornstein-Uhlenbeck process), corresponding to equation \eqref{1b} with $Y=0$ and $a=0$, is strong Feller (see Section 9.4.1 of \cite{redbook}). Therefore, if Hypothesis \ref{h3b1} holds then the solution $X$ to \eqref{1b} is strong Feller as well by Lemma 20.1 in \cite{PrZ1}.
\par
We note also that by Hypothesis \ref{h3b1} and the strong continuity of the semigroup $\la S_t\ra$ the set 
\[Q_t^{1/2}(H)=\bigcup_{s\in(0,t]}Q_s^{1/2}(H)\]
is dense in $H$ for every $t>0$. Hence, $N\la S_tx,Q_t\ra={\mathcal Law}\la X^{OU,x}_t\ra$ has a support on the whole of $H$ for every $t>0$. Therefore $X^x_t$ is irreducible for every $t>0$ and $x\in H$ by Lemma 20.1 in \cite{PrZ1}. 
\par
Thus we obtain the following immediate consequence of Theorem \ref{th1}.

\begin{corollary}\label{cor1g}
Assume Hypotheses \ref{h1}, \ref{h2} and \ref{h3b1} and moreover assume that the semigroup $\la S_t\ra$ is exponentially stable. Then equation \eqref{1b} has a unique invariant measure $\mu$ that satisfies \eqref{eq3}. 
\end{corollary}
In particular, invoking Remark \ref{22a} we obtain 
\begin{corollary}\label{2g}
Let $U=H$ and $\E\left|Y_1\right|^p<\infty$ for certain $p>0$. Moreover, assume \eqref{eqh2} and Hypothesis \ref{h3b1} and that the semigroup $\la S_t\ra$ is exponentially stable. Then equation \eqref{1b} has a unique invariant measure $\mu$ that satisfies \eqref{eq3}. 
\end{corollary}
\subsection{Heat equation with L\'evy noise with singular Gaussian component}\label{singular}
In this section we will apply Theorem \ref{th1} to the stochastic heat equation 
\begin{equation}\label{heat1}
\left\{\begin{array}{ll}
dX_t(\xi)=\la\frac{\partial^2 X_t}{\partial\xi^2}(\xi)+a(\xi)\ra dt+Q^{1/2}dW_t(\xi)+dY_t(\xi),&\xi\in(0,\pi),\,\, t>0\\
&\\
X_0(\xi)=x(\xi),&\xi\in[0,\pi]\\
&\\
X_t(0)=X_t(\pi)=0,&t>0,
\end{array}\right.
\end{equation}
with an unbounded operator $Q=\la-\frac{\partial^2}{\partial\xi^2}\ra^{\delta/2}$, vector $a\in U$, and a cylindrical $\alpha$-stable process $Y$. It is well known that equation \eqref{heat1} can be rewritten as a version of an abstract equation \eqref{1b}. Indeed, let $H=L^2(0,\pi)$ and 
\[A=\frac{\partial^2}{\partial\xi^2},\quad \mathrm{dom}(A)=H^2(0,\pi)\cap H_0^1(0,\pi).\]
Then $A$ has eigenfunctions 
\[e_k(\xi)=\sqrt{\frac{2}{\pi}}\sin k\xi,\quad k\ge 1,\]
and eigenvalues 
\[-\lambda_k=-k^2,\quad k\ge 1.\]
We assume that for a certain $\delta\in\R$
\begin{equation}\label{311}
Q=\sum_{k=1}^\infty q_ke_k\otimes e_k,\quad q_k=k^\delta.
\end{equation}
Then 
\[\mathrm{for\,\, every\,\,} t>0,\quad \sup_{k}\sqrt{\frac{\lambda_k}{q_k}}e^{-\lambda_k t}<\infty,\]
and therefore, by Proposition 9.30 in \cite{redbook} Hypothesis \ref{h3b1} holds. Note that \eqref{eqh2} holds for $A$ and $Q$ defined above iff 
\begin{equation}\label{312}
\delta<1.
\end{equation} 
Indeed, \eqref{312} is satisfied iff 
\[\sum_{k=1}^\infty\frac{q_k}{\lambda_k}<\infty,\]
which is equivalent to \eqref{eqh2}, comp. Proposition 9.30 in \cite{redbook}. 
In particular, for $\delta\in(0,1)$ the operator $Q=(-A)^{\delta/2}$ in \eqref{311} is unbounded but condition \eqref{eqh2} holds.\\
Let 
\begin{equation}\label{313}
Y(t)=\sum_{k=1}^\infty b_kY_k(t)e_k,
\end{equation}
where $\la Y_k\ra$ is a sequence of real-valued, independent, symmetric and $\alpha$-stable processes with $\alpha\in(0,2)$. We will assume that the sequence of non-negative real numbers $\la b_k\ra$ enjoys the following property.\\
The set 
\[\mathcal K=\left\{k\in\mathbb N:\,b_k>0\right\},\]
is infinite and there exists $\gamma\in\R$ such that 
\begin{equation}\label{314}
b_k=k^\gamma,\quad k\in\mathcal K.
\end{equation}
Observe, that if 
\begin{equation}\label{315}
\gamma<\frac{1}{\alpha},
\end{equation}
then Condition \eqref{53} holds and by results in Section \ref{stable} Hypothesis \ref{h2} is satisfied for $p<\alpha$. In particular, if $0\le \gamma<\frac{1}{\alpha}$ then the process $Y$ defined by \eqref{313} and \eqref{314} is a cylindrical noise taking values in the space $U_{\frac{1}{\alpha}}$ (see Section \ref{stable}). Note that if $\mathcal K\neq\mathbb N$ then the process $Y$ is degenerate. Let 
\[a=\sum_{k=1}^\infty a_ke_k\in U_{\frac{1}{2}}\, ,\]
that is
\begin{equation}\label{316}
\sum_{k=1}^\infty\frac{a_k^2}{k^2}<\infty\, .
\end{equation}
Then for every $t\ge 0$ 
\[\int_0^t\left|S_sa\right|^2ds<\infty\, ,\]
and \eqref{eqh1} follows. Finally, as a consequence of Corollary \ref{cor1g} we obtain the following 
\begin{corollary}\label{38}
Under assumptions \eqref{311} - \eqref{316} equation \eqref{heat1} has a unique stationary distribution $\mu$ and for any $p\in(0,\alpha)$ exponential estimate \eqref{35} holds. 
\end{corollary}
\section{Proofs}\label{proofs}
\subsection{Proof of Remark \ref{22a}}
If $p>1$, then the process $\la Y_t - \mathbb E Y_t\ra_{t\geq 0}$ is an $\la \mathcal  F_t\ra$-martingale and \eqref{eqh3}  and \eqref{eqh4} follow easily.\\
Let $0 < p \leq 1$. It has been proved in \cite{kruglov}, see also \cite{jurek} for a more general result, that $\mathbb E \left|Y_1\right|^p <\infty$ iff
\begin{equation}\label{33}
\int_{|x|\geq 1} |x|^p M(dx) < \infty,
\end{equation}
where $M$ is the L\'evy measure of $Y$. 
We can represent
$$Y_t = Y_t^1 +Y_t^2,$$
where $Y^1$ is a L\'evy process corresponding to the L\'evy measure $M$ restricted to the ball $B_1$ and $Y^2$ is a L\'evy process corresponding to $M^{(2)}$, the L\'evy measure $M$ restricted to $H\setminus B_1$, and moreover $Y^1$ and $Y^2$ are independent. Therefore \eqref{eqh3} is satisfied iff \eqref{eqh3} holds for $Y^2$. By \cite{PZ} we have
$$\int_0^t S_s dY^2_s =\int_0^t \int_{H\setminus B_1} S_s x N(ds, dx),$$
where $N$ is the Poisson measure corresponding to $Y^2$, in particular $N$ takes nonnegative integer values and hence $\mathbb E |N(B)|^p \leq \mathbb E |N(B)| = M^{(2)}(B)$ for $B\in\mathcal  B(H)$. Therefore as in \cite{BZ} we obtain the estimate:
\[\begin{aligned}
\mathbb E\left|\int_0^t S_s d Y_s^2\right|^p &\leq \int_0^t \int_{H\setminus B_1} \left|S_s x\right|^p M^{(2)} (dx) ds\\
&\leq \la\int_0^t\left\|S_s\right\|^p ds\ra\la \int_{H\setminus B_1} |x|^p M(dx)\ra < +\infty,\end{aligned}\]
where the last inequality follows from \eqref{33} and \eqref{eqh3} and \eqref{eqh4} hold.
\subsection{Proof of Proposition \ref{pr1}}
We will consider only the case of $p\in(0,1)$ which is less standard. The proof for $p\ge 1$ is similar. We recall that for $p\in(0,1)$ the space $L^p(\Omega, H)$ is a linear complete metric space with the distance 
\[d_p(\xi,\eta)=\E|\xi-\eta|^p.\]
\emph{Proof of (1)}. For $u>0$, $v\ge 0$ we have (e.g. \cite{Ch})
\begin{equation}\label{eq7}
\int_v^{u+v}S_{u+v-s}dY_s\stackrel{D}{=}\int_0^uS_{u-s}dY_s\stackrel{D}{=}\int_0^uS_sdY_s,
\end{equation}
where $\,\stackrel{D}{=}\,$ stands for the equality of probability distributions. Hence, invoking the fact that for certain constants $C_A\ge 1$ and $\delta$
\[\left\|S_t\right\|\le C_Ae^{\delta t},\quad t\ge 0,\]
we find that for a fixed $t>0$ 
\begin{equation}\label{eq8}
\begin{aligned}
\E\left|Y_A(t)\right|^p&=\E\left|\int_0^tS_sdY_s\right|^p\\
&\le \sum_{k=0}^{[t]-1}\E\left|\int_k^{k+1}S_sdY_s\right|^p+\left\|S_{[t]}\right\|^p\E\left|\int_0^{t-[t]}S_sdY_s\right|^p\\
&\le \gamma_pC_A^p\sum_{k=0}^{[t]}e^{p\delta k}<\infty,
\end{aligned}
\end{equation}
and \eqref{eq2} follows with 
\[C(t,p,A)=C_A^p\sum_{k=0}^{[t]}e^{p\delta k}.\]
In \eqref{eq8} estimates follow from the triangle inequality and \eqref{eq7}. \\
\emph{Proof of (2).} Fix $t\ge 0$. Then for $h>0$ 
\[Y_A(t+h)-Y_A(t)=\int_t^{t+h}S_{t+h-s}dY_s+\int_0^t\la S_{t+h-s}-S_{t-s}\ra dY_s.\]
Hence, by triangle inequality and \eqref{eq7} 
\[
\begin{aligned}
\E\left|Y_A(t+h)-Y_A(t)\right|^p&\le \E\left|\int_0^hS_sdY_s\right|^p+\E\left|S_hY_A(t)-Y_A(t)\right|^p\\
&=J_1(h)+J_2(t,h),
\end{aligned}\]
and by \eqref{eqh4} of Hypothesis \ref{h2}, 
\begin{equation}\label{eq9}
\lim_{h\downarrow 0}J_1(h)=0.
\end{equation}
Concerning $J_2(t,h)$ observe that for $\P$-a.e. $\omega$ 
\[\lim_{h\downarrow 0}S_hY_A(t,\omega)=Y_A(t,\omega),\]
by the strong continuity of the semigroup $\la S_t\ra$. Moreover, for $h\in(0,1]$ 
\[\begin{aligned}
\left|\la S_h-I\ra Y_A(t,\omega)\right|&\le\la \left\|S_h\right\|+1\ra\left|Y_A(t,\omega)\right|\\
&\le C_A\la 2+e^\delta\ra\left|Y_A(t,\omega)\right|,
\end{aligned}\] 
Therefore, by part (1) of the Proposition and the Dominated Convergence we obtain 
\begin{equation}\label{eq10}
\lim_{h\downarrow 0}J_2(t,h)=0,
\end{equation} 
and the right-continuity of $Y_A$ in $p$-th mean follows. Similarly, to prove left-continuity, fix $t>0$ and let $0<h<t\wedge 1$. Then, using \eqref{eq7} we obtain 
\[\begin{aligned}
\E\left|Y_A(t)-Y_A(t-h)\right|^p&=\E\left|\la S_h-I\ra Y_A(t-h)+\int_{t-h}^tS_{t-s}dY_s\right|^p\\
&\le J_1(h)+\E\left|\la S_h-I\ra\int_0^{t-h}S_sdY_s\right|^p\\
&=J_1(h)+J_3(t,h),
\end{aligned}\]
and
\[\begin{aligned}
J_3(t,h)&\le\E\left|\la S_h-I\ra\int_0^tS_sdY_s\right|^p+\left\|S_h-I\right\|^p\E\left|\int_{t-h}^tS_sdY_s\right|^p\\
&\le J_2(t,h)+C_A^p\la 2+e^\delta\ra^pJ_1(h). 
\end{aligned}\]
Therefore, by \eqref{eq9} and \eqref{eq10}
\[\lim_{h\downarrow 0}J_3(t,h)=0\]
and part (2) of  Proposition \ref{pr1} follows. \\
\emph{Proof of (3). } The Chebyshev inequality and part (2) of Proposition \ref{pr1} yield the stochastic continuity of the process $Y_A$. Since $Y_A$ is $\la \mathcal F_t\ra$-adapted, part (3) follows from Proposition 3.6 of \cite{redbook}.
\subsection{Proof of Proposition \ref{pr2}}
We follow the proof of Theorem 5.4 in \cite{PrZ3}. Recall that we can choose a predictable version of the process $Z_A$. Uniqueness of solutions to \eqref{1} is an easy consequence of the Gronwall Lemma. To prove existence observe first that under Hypothesis \ref{h3}, for fixed measurable function $f:[0,T]\to H$ and $x\in H$, the deterministic equation 
\[y(t)=S_tx+\int_0^tS_{t-s}F(y(s)+f(s))\,ds\]
has a unique continuous solution $y:[0,T]\to H$, hence by part (3) of Proposition \ref{pr1} there exists a continuous $H$-valued and $\la\mathcal F_t\ra$-adapted process $\la V_t\ra$ such that $\P$-a.s. 
\[V_t^x=S_tx+\int_0^tS_{t-s}F\la V_s^x+Z_A(s)\ra ds,\quad t\ge 0.\]
Consequently, the process $X_t^x=V_t^x+Z_A(t)$, $t\ge 0$, is a predictable solution to \eqref{1a}. By Proposition \ref{pr1} and Remark \ref{rem1} we have $\E\left|X_t^x\right|^p<\infty$. The Markov property is proved in Theorem 5.4 of \cite{PrZ3} and the Feller property follows from the Gronwall Lemma.  
\subsection{Proof of Theorem \ref{th1}}
\emph{Step 1. } We show first that there exist positive constants $p,\kappa, c_1,c_2$ such that 
\begin{equation}\label{eq6}
\E\left|X_t^x\right|^p\le c_1|x|^pe^{-\kappa pt}+c_2,\quad t\ge 0,\,\,x\in H.
\end{equation}
To this end we note that by Hypothesis \ref{h3}
\[C_F=\sup_{x\in H}|F(x)|<\infty,\]
and since the semigroup $\la S_t\ra$ is exponentially stable, there exist constants $C_A,\kappa>0$ such that 
\[\left\|S_t\right\|\le C_A e^{-\kappa t},\quad t\ge 0.\]
Therefore, using Remark \ref{rem1} and \eqref{eq8} for $Z_A(t)$ we obtain 
\[\begin{aligned}
\E\left|X_t^x\right|^p&\le \left|S_tx\right|^p+\E\la\int_0^t\left\|S_{t-s}\right\|\cdot\left|F\la X_s^x\ra\right| ds\ra^p+\tilde{\gamma_p}C_A^p\sum_{j=0}^{[t]}e^{-\kappa pj}\\
&\le C_A^pe^{-\kappa pt}|x|^p+C_A^pC_F^p\kappa^{-p}+\tilde{\gamma_p}C_A^p\la 1-e^{-\kappa p}\ra^{-1},
\end{aligned}\]
and \eqref{eq6} follows.\\
\emph{Step 2. } In order to prove the existence of an invariant measure for equation \eqref{1} we will use the method of \cite{DZ1} and \cite{PXZ} to show that the family $\left\{\mu_t^0:\, t\ge 1\right\}$ is tight. Indeed, for $t\ge 1$ 
\[\begin{aligned}
X_t^0&=S_1X^0_{t-1}+\int_{t-1}^t S_{t-s}F\la X_s^0\ra ds+\int_{t-1}^t S_{t-s}dZ_s\\
&=X_{1,t}+X_{2,t}+X_{3,t}.
\end{aligned}\]
Fix arbitrary $\epsilon\in(0,1)$. Since for any $t\ge 1$ 
\[X_{3,t}\stackrel{D}{=}\int_0^1 S_sdZ_s,\]
there exists a compact set $K_3=K_3(\epsilon)\subset H$ such that 
\[\P\la X_{3,t}\in K_3\ra\ge 1-\frac{\epsilon}{2},\quad\mathrm{for\,\, all}\,\, t\ge 1.\]
Next, by \eqref{eq6} 
\[\sup_{t\ge 0}\E\left|X_t^0\right|^p\le c_2<\infty,\]
hence by the Chebyshev inequality there exists $r=r_\epsilon>0$ such that 
\[\P\la X_t^0\in B_r\ra\ge 1-\frac{\epsilon}{2},\quad\mathrm{for\,\, all}\,\, t\ge 0,\]
where $B_r$ stands for the closed ball $B_r=\{x\in H:\, |x|\le r\}$. 
Since $S_1:H\to H$ is compact, the set
\[K_1=K_1(\epsilon)=S_1\la B_r\ra,\]
is compact in $H$ and 
\[\P\la S_1X_{t-1}^0\in K_1\ra\ge \P\la X_{t-1}^0\in B_r\ra\ge 1-\frac{\epsilon}{2},\]
for all $t\ge 1$. Finally, since $\la S_t\ra$ is a compact semigroup, by Lemma 6.1.4 in \cite{DZ1} the operator $G:L^2([0,1];H)\to H$ defined as
\[Gf=\int_0^1S_{1-s}f(s)\,ds\]
is compact. Then, invoking the boundedness of $F$ we find that there exists a compact set $K_2$ such that 
\[X_{2,t}\in K_2,\quad\mathrm{for\,\, all\,\,} t\ge 1.\]
Let 
\[\tilde K=K_1+K_2+K_3=\left\{x_1+x_2+x_3:\, x_i\in K_i\right\}.\]
Clearly, $\tilde K$ is a compact subset of $H$ and 
\[\P\la X_t^0\notin \tilde K\ra\le\sum_{i=1}^3\P\la X_{i,t}\notin K_i\ra\le\epsilon,\quad t\ge 1.\]
Therefore, the set $\left\{\mu_t^0;\, t\ge 1\right\}$ is tight and an invariant measure for \eqref{1} exists by the Krylov-Bogoliubov Theorem (e.g. Theorem 16.2 of \cite{PZ}).\\
\emph{Step 3. }To show \eqref{eq3} we will use the following Theorem 12.1 proved in \cite{GM}. 
\begin{theorem}\label{gm}
If an $H$-valued Markov process $\la X_t^x\ra$ is strongly Feller and irreducible, admits an invariant measure $\mu$, satisfies \eqref{eq6} and 
\begin{equation}\label{eq9a}
\la\forall_{r>0}\exists_{T>0}\exists_{\mathrm{compact\,}K\subset H}\ra\inf_{x\in B_r}\mu_T^x(K)>0,
\end{equation}
then \eqref{eq3} holds. 
\end{theorem}
Since the process $\la X_t^x\ra$ is strongly Feller and irreducible, it remains to show \eqref{eq9a}. 
Fix $r>0$ and $T>0$. Then there exists a compact $K\subset H$ such that $\mu_T^0(K)>\frac{1}{2}$. By the strong Feller property all the measures in the family $\left\{\mu_T^x;\, x\in H\right\}$ are equivalent, hence 
\begin{equation}\label{eq10a}
\mu_T^x(K)>0,\quad x\in H. 
\end{equation}
To prove \eqref{eq9a} we use a nice observation of  \cite{MS}: since $B_r$ is compact in the weak topology $\tau_w$ of $H$, \eqref{eq9a} will follow from \eqref{eq10a} and the fact that the function
\begin{equation}\label{eq11}
\la B_r,\tau_w\ra\ni x\longrightarrow \mu_T^x(K)\in\R\quad\mathrm{is\,\, continuous}.
\end{equation}
It remains to prove \eqref{eq11}. To this end, take a sequence $\la x_n\ra\subset B_r$, such that $x_n\stackrel{w}{\longrightarrow}x$. Since $S_t$ is compact for any $t>0$ we have 
\[\lim_{n\to\infty}\left|S_t\la x_n-x\ra\right|=0.\]
Next, by the Lipschitz continuity of $F$ we obtain 
\[\left|X_t^{x_n}-X_t^x\right|\le\left|S_t\la x_n-x\ra\right|+L\int_0^t\left\|S_{t-s}\right\|\cdot\left|X_s^{x_n}-X_s^x\right|\, ds,\]
hence for a certain constant $L_1>0$ the Gronwall Lemma yields 
\[\left|X_t^{x_n}-X_t^x\right|\le\left|S_t\la x_n-x\ra\right|+L_1\int_0^t\left|S_s\la x_n-x\ra\right|\,ds,\]
and thereby, using Dominated Convergence we obtain 
\[\lim_{n\to\infty}\left|X_t^{x_n}-X_t^x\right|=0,\quad\mathrm{uniformly\,\, in\,\,}\omega.\]
Now, by the Dominated Convergence 
\begin{equation}\label{eq12}
P_tf\la X_t^{x_n}\ra=\E f\la X_t^{x_n}\ra\longrightarrow\E f\la X_t^x\ra=P_tf(x)
\end{equation}
for any $t>0$ and any $f\in C_b(H)$. Putting in \eqref{eq12} $t=\frac{T}{2}$, and $f=P_tI_K\in C_b(H)$ by the strong Feller property, we find that 
\[\mu_T^{x_n}(K)=P_TI_K\la x_n\ra\longrightarrow P_TI_K(x)=\mu_T^x(K),\]
which implies \eqref{eq11} since $\la B_r,\tau_w\ra$ is metrisable.

\end{document}